\documentclass[a4paper,12pt]{article}

\usepackage[dvips]{graphics}
\usepackage{epsf}
\usepackage{graphicx}
\usepackage{epsfig}
\usepackage{amsmath}
\usepackage[psamsfonts]{amssymb}
\usepackage[psamsfonts]{eucal}
\usepackage{amsthm}
\usepackage{textcomp}
\usepackage{mathptmx}
\usepackage[scaled]{helvet}

\theoremstyle{plain}
\newtheorem{thm}[subsection]{Theorem}

\newtheorem{cor}[subsection]{Corollary}
\newtheorem{prop}[subsection]{Proposition}
\theoremstyle{definition}
\newtheorem{defn}[subsection]{Definition}

\newtheorem{rmk}[subsection]{Remark}
\newtheorem{rmks}[subsection]{Remarks}
\newtheorem{obs}[subsection]{Observation}

\def\bb#1#2{\left\{#1,#2\right\}}
\title {Poisson quasi-Nijenhuis structures with background}
\author{\large{\textsc{Paulo Antunes}}\\
  \small{\textsc{CMUC, Department of Mathematics}}\\
  \small{\textsc{University of Coimbra}}\\
  \small{\textsc{3001-454 Coimbra}}\\
  \small{\textsc{Portugal}}\\
  \footnotesize{pantunes@mat.uc.pt}}
\date{}
\begin{document}
\maketitle











\begin{abstract}
We define the Poisson quasi-Nijenhuis structures with background on Lie algebroids
and we prove that to any generalized complex
structure on a Courant algebroid which is the double of a Lie algebroid
is associated  such a structure.
We prove that any Lie algebroid
with a Poisson quasi-Nijenhuis structure with background constitutes,
with its dual, a quasi-Lie bialgebroid.
We also prove
that any pair $(\pi,\omega)$ of a Poisson bivector and a 2-form
induces a Poisson quasi-Nijenhuis structure with background and we
observe that particular cases correspond to already known
compatibilities between $\pi$ and $\omega$.

\end{abstract}



\section*{Introduction}

The aim of this work\footnote{This paper was presented as a poster at the ``Poisson 2008'' conference at the EPFL in Lausanne in July 2008.} is to define the notion of Poisson quasi-Nijenhuis structure on a Lie algebroid with a (closed) $3$-form background. The Poisson quasi-Nijenhuis structures (without background) were introduced by Sti\'{e}non and Xu in~\cite{SX07} on the tangent Lie algebroid and then on any Lie algebroid by Caseiro et al. in ~\cite{CdNNdC08}. The case with background was already studied by Zucchini~\cite{Zuc07} but we remarked that a condition is missing in the definition proposed there. This extra condition was already in~\cite{SX07} and appears naturally in this work when we ask for some structures to be integrable (or some brackets to verify the Jacobi identity).

\

In this work we will use a supermanifold approach~\cite{Vai97,Roy022}
to describe Lie algebroid structures. Let us consider a vector bundle $A \to M$
and change the parity of the fiber coordinates
(considering them odd), then we obtain a supermanifold denoted by $\Pi
A$. The algebra of functions on $\Pi A$, which are polynomial in the
fibre coordinates, is denoted by $C^\infty(\Pi A)$ and coincides with
$\Omega(A):= \Gamma(\bigwedge^{\bullet}A^*)$, the exterior algebra of
$A$-forms. Let consider a Lie algebroid structure on $A$ given by $d$,
a degree $1$ derivation of $\Omega(A)$ such that $d^2=0$. In this
supermanifold setting, $d$ is a vector field on $\Pi A$ and can be
seen as the derivation defined by a hamiltonian on $\Pi A$, i.e. an element $\mu \in
C^\infty(T^*\Pi A)$. Then $d = \{\mu, . \})$, where the so-called big bracket~\cite{K-S92},
$\bb ..$, is the canonical Poisson bracket on the symplectic supermanifold $T^*\Pi A$.
The condition $d^2=0$ is equivalent to $\bb \mu\mu=0$.

To each $f\in C^\infty(T^*\Pi A)$ is associated a bidegree
$(\epsilon,\delta)$. In fact, since using Legendre transform
(see~\cite{Roy99}) $T^{*}(\Pi A)\cong T^{*}(\Pi A^{*})$, we can define
$\epsilon$ (resp. $\delta$) as the polynomial degree of $f$ in the
fibre coordinates of the vector bundle $T^{*}(\Pi A)\to \Pi A$
(resp. $T^{*}(\Pi A)\to \Pi A^*$). We define the shifted bidegree of
$f$ as the pair $(\epsilon -1,\delta -1)$ and the total shifted
bidegree as the sum $(\epsilon -1)+(\delta -1)=\epsilon +\delta
-2$. Then a Lie algebroid structure in $A$ is a hamiltonian $\mu \in
C^\infty(T^*\Pi A)$ of shifted bidegree $(0,1)$ such that $\bb \mu\mu
=0$.

\

Instead of the expression ``with background" used here, some authors use ``twisted", or in Physics, ``H-flux". In this work, our choice was motivated by the result of the proposition~\ref{TPqN_of_(pi,omega)_as_twist}. In fact, we prove there that a particular class of Poisson quasi-Nijenhuis structure with background is obtained by twisting, in a way explained in~\cite{Roy02,Ter06,K-S07}, a Lie algebroid structure by a Poisson bivector and then by a $2$-form. Therefore, to avoid confusion, we will use the word ``twist'' only when we are dealing with twisting by a $2$-form or a bivector as in~\cite{Roy02,Ter06,K-S07}.

\

The content of this paper is as follows. In the first section we recall some basic definitions about Nijenhuis tensors, Poisson bivectors and Poisson Nijenhuis structures on a Lie algebroid and give the corresponding expression in the supermanifold approach. Then, in the second section, we introduce the notion of Poisson quasi-Nijenhuis structure with a $3$-form background $H$, on a Lie algebroid $(A,\mu)$. We prove that any complex structure (or more generally any c.p.s. structure, see definition~\ref{def_cps}) on $(A\oplus A^*,\mu + H)$ induces such a structure. In the third section we generalize a result from~\cite{SX07,CdNNdC08} and prove that any Poisson quasi-Nijenhuis structure with background on $A$ induces a Lie quasi-bialgebroid on $(A^*,A)$. Finally, in the last section we study Poisson quasi-Nijenhuis structures with background defined by a pair $(\pi,\omega)$ of a Poisson bivector and a $2$-form. We observe that already known compatible pairs such that
complementary $2$-forms for Poisson bivectors~\cite{Vai96}, Hitchin pairs~\cite{Cra04} and $P\Omega$-structures or $\Omega N$-structures~\cite{MM84}
 are all particular examples of Poisson quasi-Nijenhuis structures with background.

\section{Basic definitions}

In this section we will recall known structures such as
Poisson Nijenhuis structures on a Lie algebroid $A$ and give their expression in terms of big bracket and polynomial functions on $T^*\Pi A$.

Let $(A,[.,.],\rho)$ be a Lie algebroid over a smooth manifold $M$. The Lie algebroid structure, $([.,.],\rho)$, can be seen as a function $\mu \in C^\infty(T^*\Pi A)$, of shifted bidegree $(0,1)$, such that $\bb\mu\mu =0$.

Consider a $(1,1)$-tensor $N \in \Gamma(A\otimes A^*)$. The Nijenhuis torsion of $N$ is defined by
$$T_N(X,Y) = [NX,NY] - N\left([NX,Y] + [X,NY] - N[X,Y]\right).$$

In terms of big bracket and elements of $C^\infty(T^*\Pi A)$, the
Nijenhuis torsion is given by
\begin{equation}\label{def_TN}
T_N = \frac{1}{2}\left(\bb N{\bb N\mu} - \bb {N^2}\mu\right).
\end{equation}

If $T_N=0$, $N$ is said to be a Nijenhuis tensor and in this case we define a new Lie algebroid structure on $A$ given by
\begin{equation}\label{N-struct}
\left\{
  \begin{array}{ll}
    [X,Y]_N=[NX,Y] + [X,NY] - N[X,Y],&\quad X,Y\in \Gamma(A),\\
    \rho_N= \rho\circ N.&
  \end{array}
\right.
\end{equation}
In the supermanifold setting, the structure $([.,.]_N,\rho_N)$
is given by $\bb N\mu \in C^\infty(T^*\Pi A)$. We denote by $d_N$ the degree $1$ derivation of $\Omega(A)$ induced by this Lie algebroid structure. Then $$d_N=\bb{\bb{N}{\mu}}..$$

\medskip

A bivector $\pi \in \Gamma(\bigwedge^2 A)$ is said to be Poisson if
$[\pi,\pi]_{SN}=0$, where $[.,.]_{SN}$ is the Schouten-Nijenhuis
bracket naturally defined on $\Gamma(\bigwedge^{\bullet}
A)$. If $\pi$ is a Poisson bivector
we define a Lie algebroid structure on $A^*$ by setting

\begin{equation}\label{pi-struct}
\left\{
  \begin{array}{ll}
    [\alpha,\beta]_{\pi}=\mathcal{L}_{\pi^{\sharp}(\alpha)}\beta - \mathcal{L}_{\pi^{\sharp}(\beta)}\alpha - d(\pi(\alpha,\beta)),&\quad\quad \alpha,\beta\in \Gamma(A^*),\\
    \rho_{\pi}= \rho\circ {\pi^{\sharp}}.&
  \end{array}
\right.
\end{equation}

In the supermanifold setting, the structure $([.,.]_{\pi},\rho_{\pi})$
is given by $\bb {\pi}\mu \in C^\infty(T^*\Pi A)$.

\begin{defn}
A Poisson bivector $\pi$ and a Nijenhuis tensor $N$ are said to be \emph{compatibles} if

$$\left\{
  \begin{array}{ll}
    N\circ\pi^{\sharp}=\pi^{\sharp}\circ{}^tN,&\\
    C_{\pi,N}=0,&
  \end{array}
\right.$$
with
$$C_{\pi,N}=\left([.,.]_N\right)_{\pi} -\left([.,.]_{\pi}\right)_{{}^tN},$$
a $C^\infty(M)$-bilinear bracket on
$\Gamma(A^*)$. When $\pi$ and $N$ are compatible, $(A,\pi,N)$
constitutes a \emph{Poisson Nijenhuis Lie algebroid}.
\end{defn}

In the supermanifold setting, we have
$$C_{\pi, N}={\bb{\pi}{\bb{N}{\mu}}} + {\bb{N}{\bb{\pi}{\mu}}}.$$

\begin{thm}
If $(A,\pi,N)$ is a Poisson Nijenhuis Lie algebroid, then
$(A_N,A^*_{\pi})$ is a Lie bialgebroid,
where $A_N$ and $A^*_{\pi}$
are the Lie algebroids defined respectively by (\ref{N-struct}) and
(\ref{pi-struct}).
\end{thm}
\begin{rmk}
    When $A=TM$ and $\mu$ is the standard Lie algebroid structure, the implication of the previous theorem becomes an equivalence (see~\cite{K-S96}).
\end{rmk}

The Lie bialgebroid $(A_N,A^*_{\pi})$ induces a Courant algebroid
structure in $A\oplus A^*$~\cite{LWX97,Roy99} which is given in the
supermanifold setting by
$$S=\bb {\pi}\mu + \bb N\mu = \bb {\pi + N}\mu.$$

\medskip

\section{Poisson quasi-Nijenhuis with background and generalized geometry}

Let $S$ be a Courant algebroid structure on $A\oplus A^*$, i.e., $S \in C^\infty(T^*\Pi A)$ is of total shifted degree $1$ and $\bb SS =0$.
Consider also a $(1,1)$-tensor $J$ on $A\oplus A^*$, seen as a map $J: A\oplus A^* \to A\oplus A^*.$
We call $J$ \emph{orthogonal} if $$<J(\mathcal{X}), \mathcal{Y}> +
<\mathcal{X} , J(\mathcal{Y})> = 0,$$ for all $\mathcal{X},\mathcal{Y}
\in \Gamma(A\oplus A^*)$,
with $<.,.>$ defined by $<X+\alpha, Y+\beta>= \beta(X) + \alpha(Y)$ for all $X,Y \in \Gamma A$, $\alpha,\beta \in \Gamma A^*$.

As in the Lie algebroid case, we can define a new bracket $[.,.]_J$ deforming by $J$ the Courant structure on $A\oplus A^*$ by setting
$$[\mathcal{X},\mathcal{Y}]_J=[J\mathcal{X},\mathcal{Y}] + [\mathcal{X},J\mathcal{Y}] - J[\mathcal{X},\mathcal{Y}],$$
for all $\mathcal{X},\mathcal{Y}\in \Gamma(A\oplus A^*)$, where $[.,.]$ is the Dorfman bracket on $A\oplus A^*$. When $J$ is an orthogonal $(1,1)$-tensor on $A\oplus A^*$, this deformed bracket is given by the hamiltonian $$S_J:=\bb J\mu \in C^\infty(T^*\Pi A).$$
We define also the Nijenhuis torsion of $J$,
$$T_J(\mathcal{X},\mathcal{Y}) = [J\mathcal{X},J\mathcal{Y}] - J\left([\mathcal{X},\mathcal{Y}]_J\right),$$
for all $\mathcal{X},\mathcal{Y}\in \Gamma(A\oplus A^*)$.


\begin{prop}\label{prop_equiv_T_J/Courant}
    \begin{enumerate}
      \item The hamiltonian $S_J$ defines a Courant structure on $A\oplus A^*$ if and only if $\bb S{T_J}=0$.
      \item $J$ is a Courant morphism from $(A\oplus A^*, S_J)$ to $(A\oplus A^*, S)$ if and only if $T_J=0$.
    \end{enumerate}
\end{prop}

\begin{defn}\label{def_cps}
    An orthogonal $(1,1)$-tensor $J$, on $A\oplus A^*$, is an \emph{almost c.p.s. structure} if ${J^2 = \lambda id_{A\oplus A^*}}$, with $\lambda \in \{-1,0,1\}$. The almost c.p.s. structure $J$ is integrable when $T_J=0$.
\end{defn}
The abbreviation ``c.p.s.'' is due to Vaisman~\cite{Vai07} and corresponds to the three different structures we are considering: if $\lambda=-1$, $J$ is an almost complex structure; if $\lambda=1$, $J$ is an almost product structure; and if $\lambda=0$, $J$ is an almost subtangent structure.

As was noticed in \cite{Cra04, Vai07}, $J$ is an almost c.p.s. structure if and only if $J$ can be represented in a matrix form by
\begin{equation}\label{matrixJ}
    J\left(
        \begin{array}{c}
            X\\
            \alpha
        \end{array}
     \right)
     =\left(
        \begin{array}{cc}
            N&\pi^{\sharp}\\
            \sigma^{\flat}&-{}^tN
        \end{array}
     \right)
     \left(
        \begin{array}{c}
            X\\
            \alpha
        \end{array}
     \right)
\end{equation}
for all $X \in \Gamma(A)$ and $\alpha \in \Gamma(A^*)$, where $\pi \in \Gamma(\bigwedge^2 A), \sigma \in \Gamma(\bigwedge^2
A^*)$ and $N \in \Gamma(A \otimes A^*)$ satisfy
        $$\left\{
        \begin{array}{ll}
            &N\circ\pi^{\sharp}=\pi^{\sharp}\circ{}^tN,\\
            &\sigma^{\flat}\circ N={}^tN\circ\sigma^{\flat},\\
            &N^2 + \pi^{\sharp}\circ\sigma^{\flat}= \lambda id_A.
        \end{array}
        \right.$$

In the supermanifold setting, $$J=\pi+N+\sigma$$ in the sense that
$J(.) = \bb{.}{\pi+N+\sigma}$. Moreover, the integrability condition of an almost c.p.s. structure, $T_J=0$, is expressed by~\cite{Gra06}
\begin{equation}\label{integrability_J}
\bb{\bb JS}{J} + \lambda S=0.
\end{equation}

\medskip

Let us now consider the case where
$S=\mu+H$, where $\mu \in C^\infty(T^*\Pi A)$ defines a Lie algebroid structure on $A$, and $H\in \Gamma(\bigwedge^3 A^*)$ is a closed $3$-form. Then $\bb SS=0$ and $S$ defines a Courant algebroid structure on $A\oplus A^*$.

The goal of this section is to relate c.p.s. structures on $(A\oplus
A^*, \mu+H)$ with the Poisson quasi-Nijenhuis structures with background which we now define.


\begin{defn}\label{def_TPqN}
        A \emph{Poisson quasi-Nijenhuis structure with background} on a Lie algebroid $A$ is a quadruple $(\pi, N, \psi, H)$ where $\pi \in \Gamma(\bigwedge^2 A)$, $N \in \Gamma(A\otimes A^*)$, $\psi \in \Gamma(\bigwedge^3 A^*)$ and $H \in \Gamma(\bigwedge^3 A^*)$  are such that $N\circ\pi^{\sharp}=\pi^{\sharp}\circ{}^tN$, $d\psi=0$, $dH=0$ and verify the conditions
        \begin{equation}\label{defn-TPqN}
            \left\{
            \begin{array}{ll}
                \pi\textrm{ is a Poisson bivector},\\
                C_{\pi, N}(\alpha, \beta) = 2\,i_{\pi^{\sharp}\alpha \wedge \pi^{\sharp}\beta} H,\\
                T_N (X,Y) =  \pi^{\sharp}(i_{NX\wedge Y}H - i_{NY\wedge X}H + i_{X\wedge Y}\psi),\\
                \textrm{d}_N \psi = \textrm{d} \mathcal{H},
            \end{array}
            \right.
        \end{equation}
        for all  \mbox{$X,Y\in\Gamma(A),\ \alpha, \beta \in \Gamma(A^*)$} and where
        $\mathcal{H}$ is the $3$-form defined by
        \begin{equation}\label{def_mathcal_H}
            \mathcal{H}(X,Y,Z)= \circlearrowleft_{X,Y,Z} H(NX,NY,Z),
        \end{equation}
        for all $X,Y,Z \in\Gamma(A).$
\end{defn}

\begin{obs}
    \begin{enumerate}
      \item In terms of big bracket and elements of $C^\infty(T^*\Pi A)$, the conditions (\ref{defn-TPqN}) correspond to
        \begin{equation}\label{defn-TPqN_in_bb}
        \left\{
            \begin{array}{ll}
                \bb{\bb{\pi}{\mu}}{\pi}=0,\\
                \bb{\bb{\pi}{\mu}}{N} + \bb{\bb{N}{\mu}}{\pi} + \bb{\bb{\pi}{H}}{\pi} =0,\\
                \bb{\bb{N}{\mu}}{N} + \bb{N^2}\mu - 2\bb{\pi}{\psi} + \bb{\bb{\pi}{H}}{N} + \bb{\bb{N}{H}}{\pi} = 0,\\
                2 \bb{\bb {N}\mu}{\psi} =\bb{\mu}{\bb{N}{\bb{N}{H}} - \bb{N^2}{H}}.
            \end{array}
        \right.
        \end{equation}
      \item If $H=0$ we recover the Poisson quasi-Nijenhuis structures defined in \cite{CdNNdC08,SX07}.
      \item The last condition of (\ref{defn-TPqN}) is missing in the
        definition proposed by Zucchini \cite{Zuc07}. In our study this
        condition appears naturally and is necessary in order to
        include the case without background, described in~\cite{CdNNdC08,SX07}.
    \end{enumerate}
\end{obs}
\medskip

\begin{thm}\label{Thm_implication_cps/tPqN}
    If an endomorphism $J$, defined by (\ref{matrixJ}), is a c.p.s. structure on $(A\oplus A^*, \mu +H)$ then $(\pi, N, -d\sigma, H)$ is a Poisson quasi-Nijenhuis structure with background on $A$.
\end{thm}
\begin{proof}
    The result follows directly by writing the integrability condition
    (\ref{integrability_J}) with $J=\pi+N+\sigma$ and $S=\mu+H$. Using the
    bilinearity of $\bb ..$ and taking into account the bidegree of each
    term we obtain the following system of equations
    $$\left\{
        \begin{array}{ll}
            \bb{\bb{\pi}{\mu}}{\pi}=0,\\
            \bb{\bb{\pi}{\mu}}{N} + \bb{\bb{N}{\mu}}{\pi} + \bb{\bb{\pi}{H}}{\pi} =0,\\
            \bb{\bb{N}{\mu}}{N} + 2\bb{\pi}{\bb{\mu}{\sigma}} +\bb{\bb{\pi}{\sigma}}\mu+ \bb{\bb{\pi}{H}}{N} + \bb{\bb{N}{H}}{\pi} + \lambda \mu= 0,\\
            \bb{\bb {N}\mu}{\sigma} + \bb{\bb {\sigma}\mu}{N} + \bb{\bb{N}{H}}{N} + \bb{\bb{\pi}{\sigma}}{H} + \lambda H=0.
        \end{array}
    \right.$$

    In the last two equations of the system we now use the algebraic conditions for $J$ to be a c.p.s. structure and more precisely the condition
    $N^2 + \pi^{\sharp}\circ\sigma^{\flat}= \lambda id_A$ which is written in terms of big bracket and elements of $C^\infty(T^*\Pi A)$ as
    $$\bb{\pi}{\sigma}=N^2 -\lambda id_A.$$
    We obtain
    $$\left\{
        \begin{array}{ll}
            \bb{\bb{\pi}{\mu}}{\pi}=0,\\
            \bb{\bb{\pi}{\mu}}{N} + \bb{\bb{N}{\mu}}{\pi} + \bb{\bb{\pi}{H}}{\pi} =0,\\
            \bb{\bb{N}{\mu}}{N} + \bb{N^2}\mu + 2\bb{\pi}{\bb{\mu}{\sigma}} + \bb{\bb{\pi}{H}}{N} + \bb{\bb{N}{H}}{\pi} = 0,\\
            \bb{\bb {N}\mu}{\sigma} + \bb{\bb {\sigma}\mu}{N} + \bb{\bb{N}{H}}{N} + \bb{N^2}{H} -2 \lambda H=0.
        \end{array}
    \right.$$

    The proof is achieved after interpreting the previous system of equations as

    \begin{equation}\label{system_of_conditions_for_cps}
        \left\{
            \begin{array}{ll}
                \pi\textrm{ is a Poisson bivector},\\
                C_{\pi, N}(\alpha, \beta) = 2\,i_{\pi^{\sharp}\alpha \wedge \pi^{\sharp}\beta} H,\\
                T_N (X,Y) =  \pi^{\sharp}(i_{NX\wedge Y}H - i_{NY\wedge X}H - i_{X\wedge Y}d\sigma),\\
                2 i_N d\sigma - d(i_N \sigma)= 2 (\mathcal{H} + \lambda H),
            \end{array}
        \right.
    \end{equation}
    for all $X,Y \in\Gamma(A)$ and $\alpha, \beta \in \Gamma(A^*)$ and $\mathcal{H}$ defined by Equation (\ref{def_mathcal_H}).
\end{proof}

\begin{rmk}
    In~\cite{Vai07}, Vaisman studied the integrability of almost c.p.s. structures on $TM\oplus T^*M$ considering both the usual Courant bracket, and also the case with a $3$-form background. The conditions that he obtained in Remark 1.5 of~\cite{Vai07} coincide with the system of conditions (\ref{system_of_conditions_for_cps}).
\end{rmk}

Note that, in the previous proof, the last equation of
(\ref{system_of_conditions_for_cps}) is only a sufficient condition
for the last equation of (\ref{defn-TPqN}). We can get an equivalence
if we impose additional conditions on the quadruple $(\pi, N, \sigma, H)$.

\begin{thm}\label{Thm_equivalence_cps/tPqN}
    An endomorphism $J$, defined by (\ref{matrixJ}), is a c.p.s. structure on $(A\oplus A^*, \mu +H)$ if and only if   $(\pi, N, -d\sigma, H)$ is a Poisson quasi-Nijenhuis structure on $A$ with background $H$ such that
$$\left\{
  \begin{array}{ll}
    N^2 + \pi^{\sharp}\circ\sigma^{\flat}= \lambda id_A,\\
    \sigma^{\flat}\circ N={}^tN\circ\sigma^{\flat},\\
    2(i_N d\sigma - \mathcal{H})= d(i_N \sigma) + 2\lambda H.
  \end{array}
\right.$$
\end{thm}

\medskip

\section{Poisson quasi-Nijenhuis with background and Lie quasi-bialgebroids}

In this section we will generalize a result proved for structures
without background in \cite{SX07, CdNNdC08}. Let $(A,\mu)$ be a Lie
algebroid over a smooth manifold $M$.

\begin{defn}
A Lie quasi-bialgebroid is a triple $(A,\delta,\varphi)$ where $A$ is a Lie algebroid, $\delta$ is a degree one derivation of the Gerstenhaber algebra $\left(\Gamma(\bigwedge^{\bullet}A),\wedge,[.,.]\right)$ and $\varphi \in \Gamma(\bigwedge^3 A)$ is such that $\delta^2=[\varphi,.]$ and $\delta\varphi=0$.
\end{defn}
The main result of the section is the following

\begin{thm}
    If $(\pi, N, \psi, H)$ is a Poisson quasi-Nijenhuis structure with background on $A$ then $(A^*_{\pi},\ d_N^H,\ \psi + i_N H)$ is a Lie quasi-bialgebroid, where \mbox{$d_N^H(\alpha)=d_N(\alpha) - i_{\pi^{\sharp}(\alpha)}H,$} for all $\alpha\in\Gamma(A^*)$.
\end{thm}

\begin{proof}
    The hamiltonian on $C^\infty(T^*\Pi A)$ which induces the structure $(A^*_{\pi},\ d_N^H,\ \psi + i_N H)$ is\linebreak \mbox{$\widetilde{S}=\bb{\pi+N}{\mu+H} + \psi$}. Considering the bidegree of each term, the equation \mbox{$\bb{\widetilde{S}}{\widetilde{S}}=0$} is equivalent to

    \begin{equation}\label{aux1}
    \left\{
        \begin{array}{ll}
            \bb{\bb{\pi}{\mu}}{\bb{\pi}{\mu}}=0,\\
            \bb{\bb{\pi}{\mu}}{\bb{\pi}{H}} + \bb{\bb{\pi}{\mu}}{\bb{N}{\mu}}=0,\\
            \bb{\bb{\pi}{H}}{\bb{\pi}{H}} + \bb{\bb{N}{\mu}}{\bb{N}{\mu}} + 2\bb{\bb{\pi}{\mu}}{\bb{N}{H}} + \\
            \phantom{================}+ 2\bb{\bb{\pi}{H}}{\bb{N}{\mu}} + 2\bb{\bb{\pi}{\mu}}{\psi} = 0,\\
            \bb{\bb{\pi}{H}}{\bb{N}{H}} + \bb{\bb{N}{\mu}}{\bb{N}{H}} +\bb{\bb{\pi}{H}}{\psi} + \bb{\bb{N}{\mu}}{\psi} =0.
        \end{array}
    \right.
    \end{equation}

    It is now straightforward to observe that the system of equations (\ref{defn-TPqN_in_bb}) implies the system (\ref{aux1}).
\end{proof}

\begin{cor}
    If $(\pi, N, -d\sigma, H)$ is a Poisson quasi-Nijenhuis structure
    with background on $A$ then $\bb{\pi+N+\sigma}{\mu+H}$ is a
    structure of Lie quasi-bialgebroid on $(A^*, A)$ or equivalently
    of a Courant algebroid on $A\oplus A^*$.
\end{cor}
\begin{proof}
    If we consider $\psi=-d\sigma$ in the previous proof, then $\widetilde{S}=\bb{\pi+N+\sigma}{\mu+H}$ and we prove that $\bb {\widetilde{S}}{\widetilde{S}}=0$.
\end{proof}
\begin{rmk}
    In the corollary above, $J=\pi+N+\sigma$ is not necessarily
    integrable, i.e., the Nijenhuis torsion $T_J$ may not vanish (see
    necessary conditions in theorem
    \ref{Thm_equivalence_cps/tPqN}). But the previous corollary proves
    that the deformed structure $S_J$ defines a Courant algebroid
    structure in $A\oplus A^*$, i.e., that $\bb S{T_J}=0$ (see
    proposition \ref{prop_equiv_T_J/Courant}).
\end{rmk}

\medskip

\section{Poisson quasi-Nijenhuis with background and compatible second order tensors}

In this section we shall consider $\pi \in \Gamma(\bigwedge^2 A)$ a
Poisson bivector and a $2$-form $\omega\in \Gamma(\bigwedge^2
A^*)$. Let us denote
              $$\begin{array}{lll}
                 \pi^{\sharp}(\alpha)=\pi(\alpha,.),\quad\forall\alpha\in\Gamma(A^*),& \quad\quad & \omega^{\flat}(X)=\omega(X,.),\quad\forall X\in\Gamma(A),\\
                 N=\pi^{\sharp}\circ\omega^{\flat},& \quad\quad & \omega_N=\omega(N.,.).\\
               \end{array}$$
Then the main result of this section is the following

\begin{thm}\label{thm_TPqN induced by_(pi,omega)}
    The quadruple $(\pi, N, d\omega_N, -d\omega)$ is a Poisson quasi-Nijenhuis structure with background on $A$.
\end{thm}
\begin{proof}
    Let us denote $\psi=d\omega_N$ and $H=-d\omega$. In terms of elements of $C^\infty(T^*\Pi A)$, we have the following correspondences
    $$\left\{
        \begin{array}{ll}
            N= \bb{\omega}{\pi},\\
            \psi=\frac{1}{2}\bb{\mu}{\bb{N}{\omega}},\\
            H=\bb{\omega}{\mu}.
        \end{array}
      \right.$$
    We easily check that $\psi$ and $H$ are closed and that $N\circ\pi^{\sharp}=\pi^{\sharp}\circ{}^tN$. To prove that $(\pi, N, \psi, H)$ is a Poisson quasi-Nijenhuis structure with background we need to verify the set of conditions (\ref{defn-TPqN}) (or equivalently the conditions (\ref{defn-TPqN_in_bb})).
    \begin{enumerate}
      \item $\pi$ is a Poisson bivector by assumption.
      \item Let us start from the fact that $\pi$ is a Poisson bivector, i.e., that
        $$\bb{\bb{\pi}{\mu}}{\pi}=0$$
        and apply $\bb{\omega}{.}$ to both sides. We get
        $$\bb{\omega}{\bb{\bb{\pi}{\mu}}{\pi}}=0.$$
        We use the Jacobi identity and obtain
        $$\bb{\bb{\omega}{\bb{\pi}{\mu}}}{\pi} + \bb{\bb{\pi}{\mu}}{\bb{\omega}{\pi}} =0$$
        and using once more the Jacobi identity in the first term of l.h.s. we have
        $$\bb{\bb{\bb{\omega}{\pi}}{\mu}}{\pi} +\bb{\bb{\pi}{\bb{\omega}{\mu}}}{\pi} + \bb{\bb{\pi}{\mu}}{\bb{\omega}{\pi}} =0,$$
        which is the second condition of (\ref{defn-TPqN_in_bb})
        $$\bb{\bb{N}{\mu}}{\pi} +\bb{\bb{\pi}{H}}{\pi} + \bb{\bb{\pi}{\mu}}{N} =0.$$

      \item As above, we start from the previous condition
        $$\bb{\bb{N}{\mu}}{\pi} +\bb{\bb{\pi}{H}}{\pi} + \bb{\bb{\pi}{\mu}}{N} =0,$$
        and apply $\bb{\omega}.$ to both sides. We obtain
        $$\bb{\omega}{\bb{\bb{N}{\mu}}{\pi}} + \bb{\omega}{\bb{\bb{\pi}{H}}{\pi}} + \bb{\omega}{\bb{\bb{\pi}{\mu}}{N}} =0,$$
        and using the Jacobi identity twice we get the required equation
        $$\bb{\bb{N}{\mu}}{N} + \bb{N^2}\mu - 2\bb{\pi}{\psi} + \bb{\bb{\pi}{H}}{N} + \bb{\bb{N}{H}}{\pi} = 0.$$
      \item The way of proving this condition is the same as above. We start from the previous condition and apply $\bb{\omega}.$ to both sides. Then, using the Jacobi identity, we get
          \begin{equation}\label{eq_aux_1}
          \bb{\bb{N}{H}}{N} + \bb{N^2}{H}  - 2\bb{N}{\psi} - \bb{\bb {N^2}\omega}{\mu}= 0.
          \end{equation}
          Finally, applying $\bb{\mu}.$ we obtain
          $$\bb{\mu}{\bb{\bb{N}{H}}{N} + \bb{N^2}{H}} - 2\bb{\mu}{\bb{N}{\psi}}= 0.$$
          Using once more the Jacobi identity and the fact that $\psi$ is closed we get
          $$2 \bb{\bb {N}\mu}{\psi} = \bb{\mu}{\bb{N}{\bb{N}{H}} - \bb{N^2}{H}}.$$
    \end{enumerate}
\end{proof}

The proof of the previous theorem suggests that starting from a Poisson bivector and composing iteratively, in a certain way, with a $2$-form we get all the conditions of the definition of a Poisson quasi-Nijenhuis structure with background. The precise way to describe this fact is using the twist of a structure by a bivector or a $2$-form as in~\cite{Roy02,Ter06,K-S07}.

\begin{prop}\label{TPqN_of_(pi,omega)_as_twist}
    If we denote by $S$ the Lie quasi-bialgebroid structure induced by the Poisson quasi-Nijenhuis structure with background $(\pi, N, d\omega_N, -d\omega)$, then $S=e^{-\omega}\circ(e^{-\pi}\mu -\mu),$ or equivalently
    $$S=e^{-\omega}(\mu_\pi ),$$
    where $\mu_\pi$ is the Lie algebroid structure defined by (\ref{pi-struct}).
\end{prop}

In the next proposition we will see that the Poisson quasi-Nijenhuis structure with background $(\pi, N, d\omega_N, -d\omega)$ is induced (as shown in theorem~\ref{Thm_implication_cps/tPqN}) by a subtangent structure.

\begin{prop}
    The $(1,1)$-tensor
    $J=\left(
        \begin{array}{cc}
            N&\pi^{\sharp}\\
            {-\omega_N}^{\flat}&-{}^tN
        \end{array}
        \right)$
    is a subtangent structure (i.e., a c.p.s. structure with $\lambda=0$) on $(A\oplus A^*, \mu-d\omega)$.
\end{prop}
\begin{proof}
    Using the theorems \ref{thm_TPqN induced by_(pi,omega)} and  \ref{Thm_equivalence_cps/tPqN} we only need to prove
    $$\left\{
        \begin{array}{ll}
            N^2 - \pi^{\sharp}\circ{\omega_N}^{\flat}= 0,\\
            {\omega_N}^{\flat}\circ N={}^tN\circ{\omega_N}^{\flat},\\
            2(i_N\,d{\omega_N} + \mathcal{H})= d(i_N\,{\omega_N}).
        \end{array}
      \right.$$
    But the verification of the two first conditions is straightforward and, using the fact that $i_N\,\omega_N=i_{N^2}\omega$, the last condition is equivalent to (\ref{eq_aux_1}).
\end{proof}


In the remaining part of this section, we will see that, in the theorem \ref{thm_TPqN induced by_(pi,omega)}, if we impose
restrictions on the $2$-form $\omega$ we get already known structures
stronger than Poisson quasi-Nijenhuis with background. We also notice
that the pairs $(\pi,\omega), (\pi,N)$ and $(\omega,N)$ thus obtained
correspond to (or slightly generalize) already known compatible pairs.


\begin{cor}
    \emph{\textbf{[Poisson Nijenhuis]}}\ If $\pi \in \Gamma(\bigwedge^2 A)$ is a Poisson bivector and $\omega\in \Gamma(\bigwedge^2 A^*)$ is a closed $2$-form such that $d\omega_N=0$, then $(\pi,N)$ is a Poisson Nijenhuis structure on $A$.
\end{cor}
\begin{rmks}
    \begin{enumerate}
      \item A pair $(\pi,\omega)$ in the conditions of the corollary above is exactly what is called a $P\Omega$-structure in \cite{MM84}.
      \item The condition $d\omega_N=0$ is the compatibility condition for $(\omega,N)$ to be a Hitchin pair as it is defined in \cite{Cra04} for $A=TM$. The pair $(\omega, N)$ above is more general because $\omega$ is not necessarily symplectic.
      \item Using the fact that $\omega$ is a closed form, we can
        prove that the compatibility condition $d\omega_N=0$ is equivalent to two other known compatibility conditions:
            \begin{itemize}
              \item $\omega$ is a complementary $2$-form for $\pi$ as
                in \cite{Vai96};
              \item $(\omega, N)$ is a $\Omega N$-structure as in \cite{MM84}.
            \end{itemize}
    \end{enumerate}
\end{rmks}
Let us justify briefly the last remark. In~\cite{Vai96}, Vaisman defines $\omega$ as a complementary $2$-form for $\pi$ when
$$[\omega,\omega]_{\pi}=0,$$
where $[.,.]_{\pi}$ is the natural extension to $\Gamma(\bigwedge^{\bullet}A^*)$ of the bracket $[.,.]_{\pi}$ defined in (\ref{pi-struct}). But in terms of big bracket and elements of $C^\infty(T^*\Pi A)$, we have
$$[\omega,\omega]_{\pi}=\bb{\bb{\omega}{\bb{\pi}{\mu}}}{\omega},$$
and using the Jacobi identity twice we obtain
$$[\omega,\omega]_{\pi}=2\bb{N}{\bb{\mu}{\omega}} - \bb{\mu}{\bb{N}{\omega}},$$
which corresponds to
\begin{equation}\label{eq_aux_2}
    [\omega,\omega]_{\pi}= 2 i_N d\omega -2 d(\omega_N).
\end{equation}
In~\cite{MM84}, Magri and Morosi define a pair $(\omega, N)$ to be a $\Omega N$-structure if a particular $3$-form $S(\omega,N)$ vanishes. But we can write
\begin{equation}\label{eq_aux_3}
    S(\omega,N)=-i_N d\omega + d(\omega_N).
\end{equation}
Therefore, using (\ref{eq_aux_2}) and (\ref{eq_aux_3}) the vanishing of $d\omega_N$ is equivalent, when $d\omega=0$, to the vanishing of $[\omega,\omega]_{\pi}$ or the vanishing of $S(\omega,N)$.

\begin{cor}
    \emph{\textbf{[Poisson quasi-Nijenhuis]}}\ If $\pi \in \Gamma(\bigwedge^2 A)$ is a Poisson bivector and $\omega\in \Gamma(\bigwedge^2 A^*)$ is a closed $2$-form then $(\pi,N, d\omega_N)$ is a Poisson quasi-Nijenhuis structure on $A$ (without background).
\end{cor}

We can also define  a Poisson Nijenhuis structure with background $(\pi, N, H)$ by considering $\psi=0$ in the definition~\ref{def_TPqN}. Up to our knowledge, this structure was never studied before. We have the following result.

\begin{cor}
    \emph{\textbf{[Poisson Nijenhuis with background]}}\ If $\pi \in
    \Gamma(\bigwedge^2 A)$ is a Poisson bivector and $\omega\in
    \Gamma(\bigwedge^2 A^*)$ is a $2$-form such that $d\omega_N=0$,
    then $(\pi,N,-d\omega)$ is a Poisson Nijenhuis structure with
    background on $A$.
\end{cor}

\begin{obs}
    In the above results, the bivector $\pi$ is a true Poisson bivector. So the last structure we obtain is different from a possible compatibility between a Poisson structure with background \cite{SW01} and a Nijenhuis tensor.
\end{obs}
\medskip

\noindent\textbf{Acknowledgements}. I would like to thank Yvette Kosmann-Schwarzbach for sug\-ges\-ting these topics and for always useful discussions about this work.



\end{document}